\documentclass[10pt]{amsart}
\usepackage{amsfonts}
\usepackage{ifthen}
\usepackage{amsthm}
\usepackage{amsmath}
\usepackage{graphicx}
\usepackage{amscd,amssymb,amsthm}
\usepackage{enumerate}
\usepackage{epstopdf}
\usepackage{mathrsfs}
\usepackage{hyperref}

\newcounter{minutes}
\divide\time by 60
\newcounter{hours}
\multiply\time by 60 \addtocounter{minutes}{-\time}

\title{Bounds for the product of modified Bessel functions}

\author[\'A. Baricz]{\'Arp\'ad Baricz}
\address{Institute of Applied Mathematics, \'Obuda University, 1034 Budapest, Hungary}
\address{Department of Economics,  Babe\c{s}-Bolyai University, Cluj-Napoca 400591, Romania}
\email{bariczocsi@yahoo.com}

\author[D.J. Ma{\v{s}}irevi{\'c}]{Dragana Jankov Ma{\v{s}}irevi{\'c}}
\address{Department of Mathematics,
University of Osijek, 31000 Osijek, Croatia}
\email{djankov@mathos.hr}

\author[S. Ponnusamy]{Saminathan Ponnusamy}
\address{Indian Statistical Institute, Chennai Centre, Society for Electronic Transactions and Security,
MGR Knowledge City, CIT Campus, Taramani, Chennai 600113, India}
\email{samy@iitm.ac.in}

\author[S. Singh]{Sanjeev Singh}
\address{Department of Mathematics,
Indian Institute of Technology Madras, Chennai 600036, India}
\email{sanjeevsinghiitm@gmail.com}

\thanks{$^{\bigstar}$The research of \'A. Baricz was supported by the J\'anos Bolyai Research Scholarship of
the Hungarian Academy of Sciences (Budapest, Hungary). The second author is on leave from  IIT Madras. The research of S. Singh was supported by the
fellowship of the University Grants Commission, India. The authors are grateful to Prof. Tibor K. Pog\'any
for the fruitful discussions during his visit to Babe\c{s}-Bolyai University of Cluj-Napoca, Romania, in September 2015.}

\newtheorem{theorem}{Theorem}

\begin{document}

\def\thefootnote{}
\footnotetext{ \texttt{File:~\jobname .tex,
          printed: \number\year-\number\month-\number\day,
          \thehours.\ifnum\theminutes<10{0}\fi\theminutes}
} \makeatletter\def\thefootnote{\@arabic\c@footnote}\makeatother



\maketitle


\begin{abstract}
In this note our aim is to present some monotonicity properties of the product of modified Bessel functions of first and second kind.
Certain bounds for the product of modified Bessel functions of first and second kind are also obtained. These bounds improve
and extend known bounds for the product of modified Bessel functions of first and second kind of order zero. A new Tur\'an type inequality is also given
for the product of modified Bessel functions, and some open problems are stated, which may be of interest for further research.
\end{abstract}

\section{\bf Properties of the product of modified Bessel functions}\label{MR}
\setcounter{equation}{0}

The product of modified Bessel functions of the first and second kind, denoted by $I_{\nu}$ and $K_{\nu},$ has found many applications in mathematical physics and engineering sciences. For more details and for a large list of references the interested reader is referred to \cite{baricz,baricz1,barpon} and to the references therein. Recently, the product $I_0K_0$ appears in \cite{bordelles_truc} to estimate the number of eigenvalues of a Schr\"odinger operator, see also \cite{truc} for more details. Motivated by the fact that the product $I_{\nu}K_{\nu}$ appears in different problems of mathematical physics, in this paper starting from the papers \cite{bordelles_truc,cantrell,truc} our aim is to present some monotonicity properties for the product of modified Bessel functions of first and second kind. Some bounds for the product of modified Bessel functions of first and second kind are also obtained. These bounds improve and extend the known bounds for the product of modified Bessel functions of first and second kind of order zero. It is important to mention here that our approach is completely different than in the papers \cite{bordelles_truc,cantrell,truc} and some of our main results are proved by using some Tur\'an type inequalities for modified Bessel functions of the first and second kind, obtained in \cite{baricz,baricz1,barpon}. Finally, we note that a new Tur\'an type inequality for the product of modified Bessel functions and some open problems are stated at the end of this section, which may be of interest for further research.

\subsection{Bounds for the product of modified Bessel functions}
In this subsection we find upper bounds for the product of modified Bessel functions of first and second kind. The key tool in our proof utilizes some known bounds for Bessel function of the first kind by Landau \cite{landau} and the integral representation for product of modified Bessel functions of first and second kind \cite[p. 253]{nist}. In the sequel, we will use the so-called Landau constants $$b_L=\sqrt[3]{2}\sup\limits_{t\in\mathbb{R}_{+}}\mathrm{Ai}(t)\ \ \ \ \mbox{and}\ \ \ \ c_L=\sup\limits_{t\in\mathbb{R}_{+}}\sqrt[3]{t}J_0(t),$$ where $\mathrm{Ai}$ stands for the Airy function,
while $J_0$ stands for the Bessel function of the first kind of zero order. The next results was motivated by the Parl's upper bound
$$I_{n+1}(x)K_m(x)<x^{n-m-1}\left((n+x)(m+1+x)\right)^{\frac{m-n}{2}},$$
$n,m\in\mathbb{N},$ $m+1\leq n,$ $x>0,$ which was used by Cantrell to obtain suitable truncation and
transient errors in the computation of the generalized Marcum $Q$-function, frequently used in radar signal processing.

\begin{theorem}\label{theorem4}
Let $\mu, \nu \in \mathbb{R}$ and $x>0$.  Then the following inequalities hold true:
\begin{enumerate}
\item[\bf a.] For $\mu>\nu>0,$
\begin{equation}\label{U_bound1}
K_\nu(x)I_\mu(x)<\frac{b_{L}}{(\mu-\nu)\sqrt[3]{\mu+\nu}}.
\end{equation}
\item[\bf b.] For $\mu\ge\nu>0,$
\begin{equation}\label{U_bound2}
K_\nu(x)I_\mu(x)\le\frac{2\pi^{\frac{3}{2}}c_L}{\sqrt{3}\Gamma\left(\frac{2}{3}\right)\Gamma\left(\frac{5}{6}\right)\sqrt[3]{2x}}
\end{equation}
\item[\bf c.] For $\mu-\nu>-\frac{1}{3},\,\mu+\nu>-1,$
\begin{equation}\label{U_bound3}
K_\nu(x)I_\mu(x)\leq\frac{\Gamma\left(\frac{2}{3}\right)\,\Gamma\left(\frac{1+3(\mu-\nu)}{6}\right)c_{L}}{2^{\frac{2}{3}}\,
\Gamma\left(\frac{5+3(\mu-\nu)}{6}\right)\sqrt[3]{2x}}.
\end{equation}
\end{enumerate}
\end{theorem}

It is important to mention here that the above bounds are not sharp for small values of $x,$ however, by using the limiting forms $I_{\mu}(x)\sim \frac{e^x}{\sqrt{2\pi x}}$ and $K_{\nu}(x)\sim \sqrt{\frac{\pi}{2x}}e^{-x}$ for large values of $x$ and $\nu,\mu$ fixed, it follows that $K_\nu(x)I_\mu(x)\sim \frac{1}{2x}$ as $x\to\infty,$ and this means that the inequalities \eqref{U_bound2} and \eqref{U_bound3} are sharp for large values of $x$.

The next theorem gives lower bounds for the product of modified Bessel functions of first and second kind. Again, the proof uses the existing bounds for modified Bessel functions of second kind \cite{gaunt} and the integral representation of product of modified Bessel function of first and second kind \cite[p. 680]{grad_ryz}. In the next theorem $L_{\nu}$ stands for the modified Struve function of the first kind.

\begin{theorem}\label{theorem5}
If $x>0$ and $\nu>0,$ then the following inequalities hold true:
\begin{equation}\label{L_bound1}
I_\nu(x)K_\nu(x)>\dfrac{\Gamma(\nu)}{2x^\nu}\left[I_\nu(2x)-L_\nu(2x)\right],
\end{equation}
\begin{equation}\label{L_bound2}
I_\nu(x)K_\nu(x)>\dfrac{1}{2\nu}-\dfrac{2x\Gamma(\nu)}{\sqrt{\pi}(1+2\nu)\Gamma\left(\nu+\frac{1}{2}\right)},
\end{equation}
\begin{equation}\label{L_bound3}
I_\nu(x)K_\nu(x)\geq\frac{1}{2\nu}-\frac{x^2}{4\nu(\nu^2-1)},
\end{equation}
provided that in the last inequality $\nu>1.$
\end{theorem}

We mention that by using the limiting forms $K_{\nu}(x)\sim \frac{1}{2}\Gamma(\nu)\left(\frac{x}{2}\right)^{-\nu}$ and $I_{\nu}(x)\sim\frac{1}{\Gamma(\nu+1)}\left(\frac{x}{2}\right)^{\nu},$ where $\nu>0$ is fixed and $x\to0,$ then we clearly get that $I_{\nu}(x)K_{\nu}(x)\sim\frac{1}{2\nu}$ as $x\to0,$ which shows that the bounds in \eqref{L_bound1}, \eqref{L_bound2} and \eqref{L_bound3} are all sharp as $x\to 0.$ However, it can be shown that they are not sharp for large values of $x.$

\subsection{Bounds used in the estimation of the eigenvalues of a Schr\"odinger operator} Our second set of results are some monotonicity properties of the product of modified Bessel functions of the first kind and second kind. The following result is motivated by the inequality $$I_0(x)K_0(x)+\ln x<1,$$ where $x\in (0,1],$ which was used recently by F. Truc to bound the number of eigenvalues of a certain Schr\"odinger operator. The above inequality was recently improved by Bordell{\`e}s and J.P. Truc \cite{bordelles_truc} by changing the number $1$ to $0.55.$ Our next result improves this bound and also extends for the general case. Here $\gamma$ stands for the Euler-Mascheroni constant.

\begin{theorem}\label{theorem1}
If $\nu\geq 0,$ then the function $x\mapsto I_{\nu}(x)K_{\nu}(x)+\ln x$ is increasing on $(0,\infty).$ Consequently, if $\nu\geq 0$ and $x\in (0,1],$ then the next inequality is valid
\begin{equation}\label{cor1_ine}
I_{\nu}(x)K_{\nu}(x)+\ln x\leq I_{\nu}(1)K_{\nu}(1).
\end{equation}
Moreover, in particular, for $x\in (0,1]$, we have
\begin{equation}\label{ine1}
\ln2-\gamma<I_0(x)K_0(x)+\ln x\leq I_0(1)K_0(1),
\end{equation}
where the constants $\ln2-\gamma$ and $I_0(1)K_0(1)\simeq0.533045{\dots}$ are the best possible.
\end{theorem}

In the case of $\nu=0,$ the next inequality \eqref{ine2} was proved recently in \cite{bordelles_truc}.

\begin{theorem}\label{theorem2}
If $x>0$ and $\nu\geq \frac{1}{2},$ then
\begin{equation}\label{ine2}
I_{\nu}(x)K_{\nu}(x)-\frac{1}{2x}<\frac{1}{16x^2}.
\end{equation}
\end{theorem}

Now, for $\nu>-1,$ let us consider the function $q_{\nu}:(0,\infty)\to(0,\infty),$ defined by $$q_\nu(x)=\displaystyle\frac{I_{\nu}(x)K_{\nu}(x)}{1+|\ln x|}.$$ Recently, F. Truc \cite{truc} proved that $$\max\limits_{x\in(0,1]} q_0(x)=\lim\limits_{x\searrow 0}q_0(x)=1.$$ This result was essential in bounding the Green function of a Sturm-Liouville operator, and motivated by the result of F. Truc our aim is to show that the maximum of $q_{\nu}(x)$ for $\nu>0$ becomes actually $I_{\nu}(1)K_{\nu}(1)$ starting from a certain value of $\nu.$ Note that since $\nu\mapsto I_{\nu}(x)K_{\nu}(x)$ is decreasing on $[0,\infty)$ (see \cite{barpon}) clearly the function $\nu\mapsto q_{\nu}(x)$ is also decreasing on $[0,\infty),$ and thus we have that $q_{\nu}(x)\leq 1$ for all $\nu\geq 0$ and $x\in(0,1].$ The next result is a significant improvement of this.

\begin{theorem}\label{theorem3}
The function $q_{\nu}$ is increasing on $(0,1]$ for $\nu\geq\frac{1}{2},$ and is decreasing on $[1,\infty)$ for $\nu>-1.$ Consequently, if $x\in(0,1]$ and $\nu\geq \frac{1}{2},$ then the next inequality is valid
\begin{equation}\label{ine3}
q_\nu(x)\leq q_\nu(1)=I_\nu(1)K_\nu(1)
\end{equation}
Moreover, for $x\geq 1$ the inequality \eqref{ine3} is valid for all $\nu>-1.$
\end{theorem}

\begin{figure}[!ht]
   \centering
       \includegraphics[width=14cm]{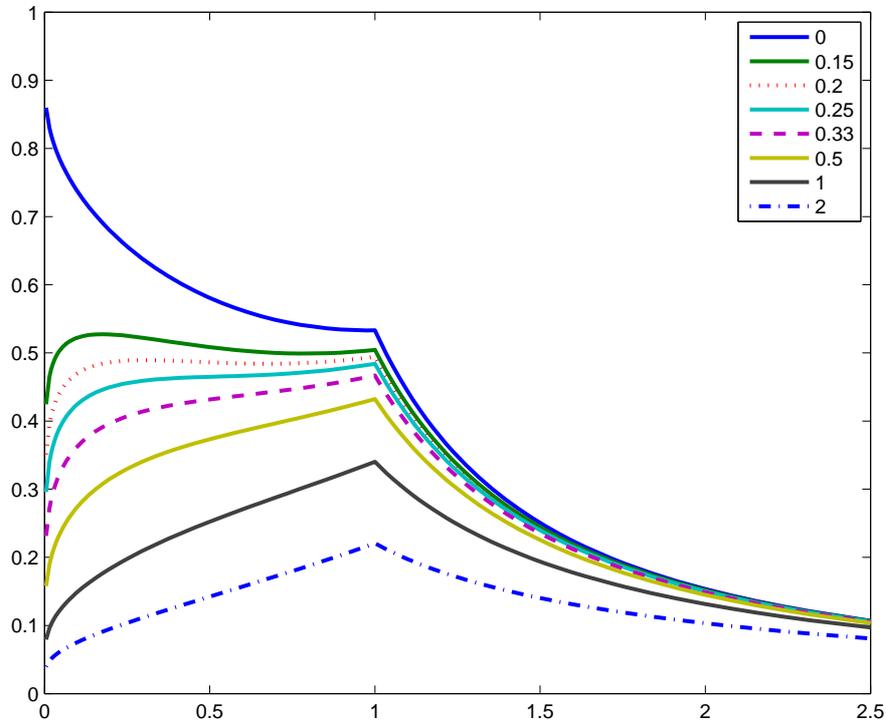}
       \caption{The graph of the function $q_{\nu}$ for $\nu\in\{0,0.15,0.2,0.25,0.33,0.5,1,2\}$ on $[0,2.5].$}
       \label{fig1}
\end{figure}

It is important to mention here that the constant $\frac{1}{2}$ is not the best possible. The next figure shows that actually there exists an $\nu^{\ast}\in(0.15,0.25)$ such that $q_{\nu}(x)\leq q_{\nu}(1)$ if and only if $\nu\geq \nu^{\ast}.$ To find this value $\nu^{\ast}$ is an open problem, which may be of interest for further research. We would like to mention that to prove the above claim a possible way it would be to study the monotonicity of $q_{\nu}$ on $(0,1].$ Figure \ref{fig1} suggest that starting from some $\nu^{\circ}\in(0.15,0.25)$ the function $q_{\nu}$ is increasing on $(0,1].$ Another open problem is to find the value of $\nu^{\circ}.$ Furthermore, we can ask whether $\nu^{\ast}$ and $\nu^{\circ}$ coincide or not.

\subsection{Tur\'an type inequalities for the product of modified Bessel functions} Our last main result complements \cite[Corollary 2]{baricz1} and it is related to the open problem posed in \cite{barpon}.

\begin{theorem}\label{th6}
The function  $\nu\mapsto\mathcal{P}_\nu(x)=\frac{\sqrt{\pi}\Gamma(\nu+\frac{1}{2})}{2x^\nu}P_\nu(x),$ where $P_\nu(x)=I_\nu(x)K_\nu(x),$ is strictly log-convex on $(-\frac{1}{2},\infty)$ for each fixed $x>0.$
\end{theorem}

As a consequence of Theorem \ref{th6} the following Tur\'an type inequality is valid for $\nu> \frac{1}{2}$ and $x>0$
\begin{equation}\label{turan}
P_\nu^2(x)-P_{\nu-1}(x)P_{\nu+1}(x)< \frac{1}{\nu+\frac{1}{2}}P_\nu^2(x).\end{equation}
We mention that since $\displaystyle1-\frac{P_{\nu-1}(x)P_{\nu+1}(x)}{P_\nu^2(x)}$ tends to zero as $x\to\infty,$ the Tur\'an type inequality \eqref{turan} is not sharp for large values of $x.$ However, for small values of $x$ such that $x^4+\left(\nu^2-\frac{1}{4}\right)x^2-\left(\nu+\frac{1}{2}\right)^2<0,$ where $\nu>\frac{1}{2},$ the inequality \eqref{turan} improves the existing Tur\'an type inequality \cite[Corollary 2]{baricz1}
$$P_\nu^2(x)-P_{\nu-1}(x)P_{\nu+1}(x)< \frac{1}{x\sqrt{x^2+\nu^2-\frac{1}{4}}}P_\nu^2(x),$$
which holds for $\nu\geq\frac{1}{2}$ and $x>0.$

Now, let us recall the following Tur\'an type inequalities (for more details see \cite{baricz1})
\begin{equation}\label{turan_I}
0<I_\nu^2(x)-I_{\nu-1}(x)I_{\nu+1}(x)<\frac{1}{\nu+1}I_\nu^2(x), \qquad \nu>-1,\ x>0
\end{equation}
and
\begin{equation}\label{turan_K}
\frac{1}{1-|\nu|} K_\nu^2(x)< K_\nu^2(x)-K_{\nu-1}(x)K_{\nu+1}(x)<0,\qquad x>0,
\end{equation}
where the left-hand side inequality is valid for all $|\nu|>1$ whereas the right-hand side inequality holds true for all $\nu\in\mathbb{R}$. By using the above inequalities \eqref{turan_I} and \eqref{turan_K} one can prove that
$$\frac{1}{1-\nu}P_\nu^2(x)<P_\nu^2(x)-P_{\nu-1}(x)P_{\nu+1}(x)<\frac{1}{\nu+1}P_\nu^2(x),\qquad x>0$$
where the left-hand side inequality is valid for all $\nu>1$ while the right-hand side inequality is true for all $\nu>0$. Observe that the right-hand side inequality in the above inequalities is stronger than the Tur\'an type inequality \eqref{turan}.

\section{\bf Proofs of the main results}\label{proof}
\setcounter{equation}{0}

\begin{proof}[\bf Proof of Theorem \ref{theorem4}]
Let us consider the integral representation of the product of modified Bessel function \cite[p. 253]{nist}
\begin{equation}\label{int1}
K_\nu(x)I_\mu(x)=\int_0^{\infty}J_{\mu+\nu}(2x\sinh t)e^{(-\mu+\nu)t}dt,\qquad \mu-\nu>-\frac{1}{2},~~\mu+\nu>-1,~~x>0.
\end{equation}
In view of the following bound for Bessel function of first kind given by Landau \cite{landau}
$$|J_\nu(x)|<b_L\nu^{-\frac{1}{3}},\qquad \nu>0,\ x\in \mathbb{R},\ b_L=0.674885{\dots},$$
it follows from \eqref{int1} that
$$
K_\nu(x)I_\mu(x)<\int_0^{\infty}b_L(\mu+\nu)^{-1/3}e^{(-\mu+\nu)t}dt=\frac{b_L}{(\mu-\nu)\sqrt[3]{\mu+\nu}},
$$
where $\mu>\nu$ and $x>0.$ This proves the inequality \eqref{U_bound1}.

Now, using another bound of Landau \cite{landau}, namely,
$$|J_\nu(x)|\leq c_L|x|^{-\frac{1}{3}} ~\mbox{ for }~\nu>0,\ x\in \mathbb{R},\ c_L=0.7857468704{\dots},$$
we obtain from \eqref{int1} that
\begin{equation}\label{U_bound}
K_\nu(x)I_\mu(x)\leq\int_0^{\infty}c_L(2x)^{-\frac{1}{3}}(\sinh t)^{-\frac{1}{3}}e^{(-\mu+\nu)t}dt=\frac{\Gamma\left(\frac{2}{3}\right)\,\Gamma\left(\frac{1+3(\mu-\nu)}{6}\right)c_{L}}{2^{\frac{2}{3}}\,
\Gamma\left(\frac{5+3(\mu-\nu)}{6}\right)\sqrt[3]{2x}},
\end{equation}
where $\mu-\nu>-\frac{1}{3},$ $\mu+\nu>-1$ and $x>0.$ This proves the inequality \eqref{U_bound3}.

If in the inequality \eqref{U_bound} we assume that $\mu\geq \nu>0,$ then we get $e^{(-\mu+\nu)t}\leq 1$ for all $t>0$ and consequently \eqref{U_bound} gives
$$K_\nu(x)I_\mu(x)\leq\int_0^{\infty}c_L(2x)^{-\frac{1}{3}}(\sinh t)^{-\frac{1}{3}}dt=\frac{2\pi^{\frac{3}{2}}c_L}{\sqrt{3}\Gamma\left(\frac{2}{3}\right)\Gamma\left(\frac{5}{6}\right)\sqrt[3]{2x}},$$
where $\mu\ge\nu>0$ and $x>0$. This proves the inequality \eqref{U_bound2}.
\end{proof}

\begin{proof}[\bf Proof of Theorem \ref{theorem5}]
To prove the inequality \eqref{L_bound1}, consider the integral representation \cite[p. 680]{grad_ryz}
\begin{equation}\label{int2}
I_\nu(x)K_\nu(x)=\frac{2x^\nu}{\sqrt{\pi}\Gamma\left(\nu+\frac{1}{2}\right)}\int_0^1t^\nu(1-t^2)^{\nu-\frac{1}{2}}K_\nu(2xt)dt,\qquad\nu>-\frac{1}{2},\ x>0.
\end{equation}
Using the lower bound for the modified Bessel function of second kind \cite[p. 383]{gaunt}
$$
K_\nu(x)>\frac{2^{\nu-1}\Gamma(\nu)e^{-x}}{x^\nu},\qquad \nu>0,\ x>0,
$$
the integral in \eqref{int2} gives
\begin{equation}\label{L_bound}
I_\nu(x)K_\nu(x)>\frac{\Gamma(\nu)}{\sqrt{\pi}\Gamma\left(\nu+\frac{1}{2}\right)}\int_0^1(1-t^2)^{\nu-\frac{1}{2}}e^{-2xt}dt
=\dfrac{\Gamma(\nu)}{2x^\nu}(I_\nu(2x)-L_\nu(2x)),
\end{equation}
where $\nu>0$ and $x>0.$ This proves the inequality \eqref{L_bound1}.

From the inequality \eqref{L_bound}, we obtain
$$I_\nu(x)K_\nu(x)>\frac{\Gamma(\nu)}{\sqrt{\pi}\Gamma\left(\nu+\frac{1}{2}\right)}\int_0^1(1-t^2)^{\nu-\frac{1}{2}}(1-2xt)dt
=\frac{\Gamma(\nu)}{\sqrt{\pi}\Gamma\left(\nu+\frac{1}{2}\right)}\left[\frac{\sqrt{\pi}\Gamma\left(\nu+\frac{1}{2}\right)}
{2\Gamma(\nu+1)}-\frac{2x}{1+2\nu}\right],$$
which is equivalent to \eqref{L_bound2}.

Now, to prove the inequality \eqref{L_bound3}, consider the inequality \cite[p. 381]{gaunt}
$$
\frac{1}{x^2}-\frac{x^{\nu-2}K_\nu(x)}{2^{\nu-1}\Gamma(\nu)}\leq \frac{1}{4(\nu-1)},\qquad x> 0,\ \nu>1,
$$
or equivalently,
$$
K_\nu(x)\geq 2^{\nu-1}\Gamma(\nu)x^{-\nu}\left[1-\frac{x^2}{4(\nu-1)}\right],\qquad x> 0,\ \nu>1.
$$
It follows from the last inequality and \eqref{int2} that
\begin{eqnarray*}
I_\nu(x)K_\nu(x)&\geq &\frac{2x^\nu}{\sqrt{\pi}\Gamma\left(\nu+\frac{1}{2}\right)}\int_0^1t^\nu(1-t^2)^{\nu-\frac{1}{2}} 2^{\nu-1}\Gamma(\nu)(2xt)^{-\nu}\left[1-\frac{(2xt)^2}{4(\nu-1)}\right]dt\\
&=&\frac{\Gamma(\nu)}{\sqrt{\pi}\Gamma\left(\nu+\frac{1}{2}\right)}\int_0^1(1-t^2)^{\nu-\frac{1}{2}} \left[1-\frac{x^2t^2}{(\nu-1)}\right]dt\\
&=&\frac{\Gamma(\nu)}{\sqrt{\pi}\Gamma\left(\nu+\frac{1}{2}\right)}\cdot\left[\frac{\sqrt{\pi}\Gamma(\nu+1/2)(2\nu^2-2-x^2)}{4(\nu-1)\Gamma(\nu+2)}\right]\\
&=&\frac{1}{2\nu}-\frac{x^2}{4\nu(\nu^2-1)},
\end{eqnarray*}
where $\nu>1$ and $x>0.$ This proves the inequality \eqref{L_bound3}.
\end{proof}

\begin{proof}[\bf Proof of Theorem \ref{theorem1}]
Let us consider the function $f_{\nu}:(0,\infty)\rightarrow \mathbb{R},$ defined by
$$f_{\nu}(x)=I_{\nu}(x)K_{\nu}(x)+\ln x.$$
Since
$$\frac{xf'_{\nu}(x)}{I_{\nu}(x)K_{\nu}(x)}=\frac{xI'_{\nu}(x)}{I_{\nu}(x)}+\frac{xK'_{\nu}(x)}{K_{\nu}(x)}+\frac{1}{I_{\nu}(x)K_{\nu}(x)},$$
in view of the Wronskian relation \cite[p. 251]{nist}
$$K_{\nu}(x)I'_{\nu}(x)-K'_{\nu}(x)I_{\nu}(x)=\frac{1}{x},$$
we get that
$$\frac{xf'_{\nu}(x)}{I_{\nu}(x)K_{\nu}(x)}=2\frac{xI'_{\nu}(x)}{I_{\nu}(x)}.$$
Now, since the function $x\mapsto \displaystyle\frac{xI'_{\nu}(x)}{I_{\nu}(x)}$ is increasing on $(0,\infty)$ for all $\nu>-1$ (see for example \cite{baricz1}) and $\displaystyle\lim_{x\rightarrow 0}\frac{xI'_{\nu}(x)}{I_{\nu}(x)}=\nu,$ we obtain for $\nu\geq0$
$$\frac{xf'_{\nu}(x)}{I_{\nu}(x)K_{\nu}(x)}=2\frac{xI'_{\nu}(x)}{I_{\nu}(x)}\geq2\nu\geq0.$$
It follows that the function $f_{\nu}$ is increasing on $(0,\infty)$ for all $\nu\geq 0.$ In particular, the function $f_0$ is increasing on $(0,1]$ and consequently on using the limits
$$
\lim_{x\rightarrow 0}f_0(x)=\ln 2-\gamma~~\mbox{and}~~
\lim_{x\rightarrow 1}f_0(x)=I_0(1)K_0(1)\approx0.533045{\dots}.
$$
the inequality \eqref{ine1} follows.
\end{proof}

\begin{proof}[\bf Proof of Theorem \ref{theorem2}]
Let us consider the function $g_{\nu}:(0,\infty)\rightarrow \mathbb{R},$ defined by
$$g_{\nu}(x)=2xI_{\nu}(x)K_{\nu}(x)-1-\frac{1}{8x}.$$
Since for $\nu\geq \frac{1}{2},$ the function $x\mapsto 2xI_{\nu}(x)K_{\nu}(x)$ is a cumulative distribution function \cite[p. 606]{hartmanwatson}, we get that $g_{\nu}$ is increasing on $(0,\infty)$ for all $\nu\geq \frac{1}{2}.$ Consequently we have that
$$g_{\nu}(x)<\lim_{x\rightarrow \infty}g_{\nu}(x)=0.$$
Thus the inequality \eqref{ine2} is indeed true for all $\nu\geq \frac{1}{2}$ and $x>0.$
\end{proof}

\begin{proof}[\bf Proof of Theorem \ref{theorem3}]
First we prove that the function $q_\nu$ is increasing on $(0,1].$ Note that $q_\nu'(x)>0$ is equivalent to
\[\left(\dfrac{xI_\nu'(x)}{I_\nu(x)}+\dfrac{xK_\nu'(x)}{K_\nu(x)}\right)(1-\ln x)+1>0.\]
In order to prove the last inequality recall that
\[\dfrac{xI_\nu'(x)}{I_\nu(x)}>\nu,\qquad \nu>-1,\ x>0\]
and take into account the next inequality (see \cite[p. 265]{laforgia})
\[\dfrac{xK_\nu'(x)}{K_\nu(x)}>-\nu-x,\qquad \nu>\frac12,\ x>0,\]
together with the fact that
$$\dfrac{xK_{\frac{1}{2}}'(x)}{K_{\frac{1}{2}}(x)}=-\frac{1}{2}-x.$$
It follows that
\[\left(\dfrac{xI_\nu'(x)}{I_\nu(x)}+\dfrac{xK_\nu'(x)}{K_\nu(x)}\right)(1-\ln x)+1>(\nu-\nu-x)(1-\ln x)+1=-x(1-\ln x)+1>0,\]
for all $x\in(0,1]$ and $\nu\geq\frac{1}{2},$ where we used the elementary inequality for the logarithm
\[\ln x>\dfrac{x-1}{x},\qquad x>0,\ x\neq 1\]
which is equivalent to
\[-x(1-\ln x)>-1,\qquad x>0,\ x\neq 1.\]

Now, we are going to prove that $q_{\nu}$ is decreasing on $[1,\infty)$ for all $\nu>-1.$ To prove this it is enough to show that
\[\left(\dfrac{xI_\nu'(x)}{I_\nu(x)}+\dfrac{xK_\nu'(x)}{K_\nu(x)}\right)(1+\ln x)-1<0,\]
where $x\geq 1$ and $\nu>-1.$ By using the Tur\'an type inequalities \cite[p. 524]{barpon}
$$\frac{xI'_{\nu}(x)}{I_{\nu}(x)}<\sqrt{x^2+\nu^2},\qquad \nu>-1,\ x>0,$$
$$\frac{xK'_{\nu}(x)}{K_{\nu}(x)}<-\sqrt{x^2+\nu^2},\qquad \nu\in \mathbb{R},\ x>0,$$
the above inequality is certainly true and consequently the inequality \eqref{ine3} is indeed true for $x\geq 1$ and $\nu>-1$.
\end{proof}

\begin{proof}[\bf Proof of Theorem \ref{th6}]
In view of the integral representation \eqref{int2} and the fact that $\nu\mapsto K_\nu(x)$ is strictly log-convex on $\mathbb{R}$ for each fixed $x > 0$ (see \cite{baricz2}), we obtain that
\begin{eqnarray*}
\mathcal{P}_{\alpha\nu_1+(1-\alpha)\nu_2}(x)&=&\int_0^1t^{\alpha\nu_1+(1-\alpha)\nu_2}(1-t^2)^{\alpha\nu_1+(1-\alpha)\nu_2-\frac{1}{2}}K_{\alpha\nu_1+(1-\alpha)\nu_2}(2xt)dt\\
&<&\int_0^1t^{\alpha\nu_1+(1-\alpha)\nu_2}(1-t^2)^{\alpha\nu_1+(1-\alpha)\nu_2-\frac{1}{2}}\left[K_{\nu_1}(2xt)\right]^{\alpha}\left[K_{\nu_2}(2xt)\right]^{1-\alpha}dt\\
&=&\int_0^1\left[t^{\nu_1}(1-t^2)^{\nu_1-\frac{1}{2}}K_{\nu_1}(2xt)\right]^{\alpha}\left[t^{\nu_2}(1-t^2)^{\nu_2-\frac{1}{2}}K_{\nu_2}(2xt)\right]^{1-\alpha}dt,
\end{eqnarray*}
which on using the H\"{o}lder-Rogers inequality
$$
\int_a^b|f(t)g(t)|dt\leq \left[\int_a^b|f(t)|^pdt\right]^{\frac{1}{p}}\left[\int_a^b|g(t)|^qdt\right]^{\frac{1}{q}},
$$
where $p > 1$, $1/p + 1/q = 1$, $f$ and $g$ are real functions defined on $[a, b]$ and
$|f|^p$, $|g|^q$ are integrable functions on $[a, b],$ gives
\begin{eqnarray*}
\mathcal{P}_{\alpha\nu_1+(1-\alpha)\nu_2}(x)&<&\left[\int_0^1t^{\nu_1}(1-t^2)^{\nu_1-\frac{1}{2}}K_{\nu_1}(2xt)dt\right]^\alpha\left[\int_0^1t^{\nu_2}(1-t^2)^{\nu_2-\frac{1}{2}}K_{\nu_2}(2xt)dt\right]^{1-\alpha}\\
&=&\left[\mathcal{P}_{\nu_1}(x)\right]^{\alpha}\left[\mathcal{P}_{\nu_2}(x)\right]^{1-\alpha},
\end{eqnarray*}
where $\alpha\in(0,1)$, $\nu_1,\nu_2\in (-\frac{1}{2},\infty)$ and $x>0.$ Hence $\nu\mapsto \mathcal{P}_\nu(x)$ is strictly log-convex on $(-\frac{1}{2},\infty)$ for each fixed $x>0$. Now, choosing $\nu_1=\nu-1$, $\nu_2=\nu+1$ and $\alpha=\frac{1}{2}$, the above inequality gives the Tur\'an type inequality
$$
\mathcal{P}_\nu^2(x)-\mathcal{P}_{\nu-1}(x)\mathcal{P}_{\nu+1}(x)<0,
$$
which is equivalent to the inequality  \eqref{turan}.
\end{proof}

\end{document}